\documentclass[12pt]{article}
\usepackage{amssymb,amsmath}
\usepackage{theorem}
\usepackage{amsfonts,pb-diagram,graphicx}
\usepackage[latin1]{inputenc}
\usepackage[numbers,sort&compress]{natbib}
\usepackage{a4}
\usepackage{tikz}
\usepackage{caption}

\usepackage{color}
\definecolor{red}{rgb}{.7,0,0} 
\definecolor{blue}{rgb}{0,0,1}

{\theorembodyfont{\slshape}
\newtheorem{theorem}{Theorem}
\newtheorem{proposition}{Proposition}

\newtheorem{lemma}[proposition]{Lemma}
}
{\theorembodyfont{\rmfamily}
\newtheorem{remark}[proposition]{Remark}

}

\def\proof{{\noindent\sc Proof. \quad}}
\newcommand{\proofof}[1]{{\noindent\sc Proof of #1. \quad}}
\def\eproof{{\mbox{}\hfill\qed}\medskip}
\newcommand\qed{{\unskip\nobreak\hfil\penalty50\hskip2em\vadjust{}
\nobreak\hfil$\Box$\parfillskip=0pt\finalhyphendemerits=0\par}}

\newcommand{\remove}[1]{}

\def\ep{\varepsilon}

\def\U_np{\U_npsilon}

\def\C{{\mathbb C}}

\def\E{\mathop{\mathbb E}}

\def\N{{\mathbb N}}

\def\P{{\mathbb P}}

\def\R{{\mathbb R}}
\def\S{{\mathbb S}}

\def\U_n{{\mathbb U}}
\def\V{{\mathcal V}}

\def\CD{{\mathcal D}}

\def\CH{{\mathcal H}}

\def\CO{{\mathcal O}}

\def\CS{{\mathcal S}}

\def\a{\alpha}
\def\NJ{{\mathrm{NJ}}}

\newcommand{\hd}{\mathcal{H}_{(d)}}

\newcommand{\czeta}{\frac{\pes{\cdot}{\zeta}^{d_i}}{\pes\zeta\zeta^{d_i}}}

\newcommand {\diag}{\mathrm{diag}}

\newcommand{\vol}{\mathrm{vol}}
\newcommand{\av}{\mathrm{av}}
\newcommand{\mufa}{\mu_{F,\hspace{0.5pt}\av}}
\newcommand{\Id}{\mathrm{Id}}
\newcommand{\Prob}{\mathop{\mathrm{Prob}}}

\def\og{\overline{g}}
\def\oU{\overline{U}}
\def\bz{{\mathbf z}}

\def\Pn{\P(\C^{n+1})}
\newcommand{\hf}{\widehat{f}}

\newcommand{\pes}[2]{\langle #1,#2\rangle}

\begin{document}
 
 \begin{title}
{Condition length and complexity for the solution of polynomial systems}
 \end{title}
\author{
Diego Armentano\thanks{Partially supported by Agencia Nacional de 
Investigaci\'on e Innovaci\'on (ANII), Uruguay, and by CSIC group 618}\\
Universidad de La Rep\'ublica\\
URUGUAY\\
{\tt diego@cmat.edu.uy}
\and
Carlos Beltr\'an\thanks{ partially suported by
the research projects MTM2010-16051 and MTM2014-57590 from 
Spanish Ministry of Science MICINN}\\
Universidad de Cantabria\\ 
SPAIN\\
{\tt beltranc@unican.es}
\and
Peter B\"urgisser\thanks{Partially funded by DFG research grant 
BU 1371/2-2}\\
Technische Universit\"at Berlin\\
%Institut f\"ur Mathematik\\
%Sekretariat MA 3-2\\
%Stra{\sz}e des 17. Juni 136
%10623 Berlin, 
GERMANY\\
{\tt pbuerg@math.tu-berlin.de}
\and
Felipe Cucker\thanks{Partially funded by
a GRF grant from the Research Grants Council of the
Hong Kong SAR (project number CityU 100813).}\\
%Department of Mathematics\\
City University of Hong Kong\\
HONG KONG\\
{\tt macucker@cityu.edu.hk}
\and
Michael Shub\\
City University of New York\\
U.S.A.\\
{\tt shub.michael@gmail.com}
}

%\date{September 30,2011}
%\subjclass{Primary 65H10, 65H20. Secondary 58C35}
\maketitle

\begin{abstract}
Smale's 17th problem asks for an algorithm which finds an approximate zero of
polynomial systems in
average polynomial time (see Smale~\cite{Smale-2000}). The main progress on Smale's problem is Beltr\'an-Pardo~\cite{Beltran_Pardo-2011} and B\"urgisser-Cucker~\cite{Buergisser_Cucker-2010}. 
In this paper we will improve on both approaches and we prove an important intermediate result. Our main results are Theorem~\ref{thm:thm1} on the complexity of a randomized algorithm which improves the result of~\cite{Beltran_Pardo-2011}, Theorem~\ref{theo:muFav} on the average of the condition number of polynomial systems which improves the estimate found in~\cite{Buergisser_Cucker-2010}, and Theorem~\ref{th:near17} on the complexity of  finding a single zero of polynomial systems. This last Theorem is the main result of~\cite{Buergisser_Cucker-2010}. We give a proof of it relying only on homotopy methods, thus removing the need for the elimination theory methods used in~\cite{Buergisser_Cucker-2010}. We build on methods developed in Armentano et al.~\cite{ABBCS}.
\end{abstract}

\section{Introduction}

Homotopy or continuation methods to solve a problem which might depend
on parameters start with a problem instance and known solution and try
to continue the solution along a path in parameter space ending at the
problem we wish to solve. We recall how this works for the solutions
of polynomial systems using a variant of Newton's method to accomplish
the continuation.
% 
% \textcolor{blue}{In this context Smale's 17th problem asks if there is an algorithm for solving polynomial systems which finds an approximate zero with probability one in average polynomial time (see Smale~\cite{Smale-2000}). The main progress on Smale's problem is Beltr\'an-Pardo\cite{Beltran_Pardo-2011} and B\"urgisser-Cucker\cite{Buergisser_Cucker-2010}. In this paper we will improve on both as well an important intermediate result. Our main results are Theorem~\ref{thm:thm1} on the complexity of a randomized algorithm which improves the result of~\cite{Beltran_Pardo-2011}, Theorem~\ref{theo:muFav} on the average of the condition number of polynomial systems which improves the estimate and simplifies the proof found in~\cite{Buergisser_Cucker-2010} and Theorem~\ref{th:near17} on the complexity of  finding a single root of polynomial systems. This last Theorem is the main result of~\cite{Buergisser_Cucker-2010}. We give a proof of it relying only on homotopy methods, thus removing the need for the elimination theory methods used in~\cite{Buergisser_Cucker-2010}. We build on methods developed in Armentano et al.\cite{ABBCS}.}

Let $\CH_{d}=\CH_{d}^{n+1}$ be the complex vector space of degree $d$
complex homogeneous polynomials in $n+1$ variables.
For $\alpha=(\alpha_0,\ldots,\alpha_n)\in\N^{n+1}$,
$\sum_{j=0}^n\alpha_j=d$, and the monomial
$z^\alpha=z_0^{\alpha_0}\cdots z_n^{\alpha_n}$, the Weyl Hermitian
structure on $\CH_{d}$ makes $\pes{z^\alpha}{z^\beta}:=0$, for
$\alpha\neq\beta$ and 
$$
     \pes{z^\alpha}{z^\alpha}:=\binom{d}{\alpha}^{-1}
    =\left(\frac{d!}{\alpha_0!\cdots \alpha_n!}\right)^{-1}.
$$

Now for $(d)=(d_1,\ldots,d_n)$ we let $\CH_{(d)}=\prod_{k=1}^n
\CH_{d_k}$. This is a complex vector space of dimension 
$$
    N:=\sum_{i=1}^n{{n+d_i}\choose{n}}.
$$
That is, $N$ is the {\em size} of a system $f\in\hd$, understood as the 
number of complex numbers needed to describe $f$.

We endow $\hd$ with the product Hermitian structure
$$
   \pes{f}{g}:=\sum_{k=1}^n \pes{f_i}{g_i},
$$ 
where $f=(f_1,\ldots,f_n)$, and $g=(g_1,\ldots, g_n)$.
This Hermitian structure is sometimes called the Weyl, Bombieri-Weyl,
or Kostlan Hermitian structure. It is invariant under unitary
substitution $f\mapsto f\circ U^{-1}$, where $U$ is a unitary
transformation of $\C^{n+1}$ (see Blum et
al.~\cite[p.~118]{Blum_Cucker_Shub_Smale-1998} for example).

On $\C^{n+1}$ we consider the usual Hermitian structure
$$
   \pes{x}{y}:=\sum_{k=0}^n x_k\,\overline{y_k}.
$$
Given $0\neq \zeta \in \C^{n+1}$, let $\zeta^\perp$ denotes
Hermitian complement of $\zeta$, 
$$
    \zeta^\perp :=\{v\in \C^{n+1}:\,\pes{v}{\zeta}=0 \}.
$$
The subspace $\zeta^\perp$ is a model for the tangent space,
$T_\zeta\mathbb{P}(\C^{n+1})$, of the projective space
$\mathbb{P}(\C^{n+1})$ at the equivalence class of $\zeta$ (which we
also denote by $\zeta$).
The space $T_\zeta\mathbb{P}(\C^{n+1})$ inherits an Hermitian 
structure from $\pes{\cdot}{\cdot}$ given by
$$
   \pes{v}{w}_\zeta :=\frac{\pes{v}{w}}{\pes \zeta\zeta}.
$$
%where $v,\,w\,\in \zeta^\perp$ represent tangent vectors at 
%$T_\zeta\mathbb{P}(\C^{n+1})$.

The group of unitary transformations $\U_n$ acts naturally on $\C^{n+1}$
by $\zeta\mapsto U^{-1}\zeta$ for $U\in\U_n$, and the Hermitian structure of
$\C^{n+1}$ is invariant under this action.

A {\em zero} of the system of equations $f$ is a point $x\in\C^{n+1}$
such that $f_i(x)=0$, $i=1,\ldots,n$.  If we think of $f$ as a mapping
$f:\C^{n+1}\to \C^n$, it is a point $x$ such that $f(x)=0$.

For a {\em generic} system (that is, for a Zariski open set of 
$f\in\CH_{(d)}$), B\'ezout's theorem states that the set of 
zeros consist of $\mathcal{D}:=\prod_{k=1}^nd_k$ complex lines
through~$0$. These $\mathcal{D}$ lines are $\mathcal{D}$ points in
projective space $\P(\C^{n+1})$. So our goal will be to approximate
one of these points, and we will use the so-called {\em homotopy} or 
{\em continuation methods}.

These methods for the solution of a system $f\in\CH_{(d)}$ proceed as
follows. 
Choose $g\in\CH_{(d)}$ and a zero $\zeta\in\P(\C^{n+1})$ of $g$.
Connect $g$ to $f$ by a path $f_t$, $0\leq t\leq 1$, in $\CH_{(d)}$
such that $f_0=g$, $f_1=f$, and try to continue $\zeta_0=\zeta$ to $\zeta_t$
such that $f_t(\zeta_t)=0$, so that $f_1(\zeta_1)=0$ (see
Beltr\'an-Shub~\cite{Beltran_Shub} for details or~\cite{Condition} for a complete discussion).

So homotopy methods numerically approximate the path
$(f_t,\zeta_t)$. One way to accomplish the approximation is via
(projective) Newton's methods. Given an approximation $x_t$
to~$\zeta_t$, define 
$$ 
   x_{t+\Delta t}:=N_{f_{t+\Delta t}}(x_t), 
$$ 
where for $h\in\CH_{(d)}$ and $y\in\Pn$ we define the {\em projective Newton's method} $N_{h}(y)$ following~\cite{Shub1993}:
$$ 
    N_{h}(y):=y-(Dh(y)|_{y^\perp})^{-1}h(y).
$$
Note that $N_h$ is defined on $\P(\C^{n+1})$ at
those points where ${Dh(y)|_{y^\perp}}$ is invertible. 

That $x_t$ is an {\em approximate zero of $f_t$ with associated (exact) zero~$\zeta_t$} 
means that the sequence of Newton iterations $N^k_{f_t}(x_t)$
converges immediately and quadratically to~$\zeta_t$.

Let us assume that $\{f_t\}_{t\in[0,1]}$ is a path in the sphere
$\S(\CH_{(d)}):=\{h\in\CH_{(d)}:\, \|h\|=1\}$. The main result of
Shub~\cite{Shub-2009}\footnote{In Shub~\cite{Shub-2009} the theorem is
actually proven in the projective space instead of the sphere, which
is sharper, but we only use the sphere version in this paper.}  
is that the $\Delta t_k$ may be chosen so that $t_0=0$, $t_k=t_{k-1}+\Delta t_k$ 
for $k=1,\ldots,K$ with $t_K=1$, such that for all $k$,  
$x_{t_k}$ is an approximate zero of $f_{t_k}$ with associated
zero $\zeta_{t_k}$, and the number $K$ of steps can be bounded as follows: 
\begin{equation}\label{eq:bezsix}
   K=K(f,g,\zeta)\leq C\,D^{3/2}\, \int_0^1 \mu(f_t,\zeta_t)\,
   \|(\dot f_t,\dot \zeta_t)\|\,dt.
\end{equation}
Here $C$ is a universal constant, $D=\max_i d_i$, 
$$
   \mu(f,\zeta) :=\begin{cases}\|f\|\,\big\|(Df(\zeta)|_{\zeta^\perp})^{-1}
   \diag(\|\zeta\|^{d_i-1}\sqrt{d_i})\big\|& \text{if $Df(\zeta)|_{\zeta^\perp}$
    is invertible }\\\infty&\text{otherwise}\end{cases}
$$
is the condition number of $f\in\CH_{(d)}$ at $\zeta\in\P(\C^{n+1})$, and
$$
   \|(\dot f_t,\dot\zeta_t)\|=(\|\dot f_t\|^2+\|\dot \zeta_t\|_{\zeta_t}^2)^{1/2}
$$
is the norm of the tangent vector to the curve in $(f_t,\zeta_t)$. The result in~\cite{Shub-2009} is not fully constructive, 
but specific constructions have been given, see~\cite{Beltran-2011} and~\cite{Dedieu_Malajovich_Shub}, 
and even programmed~\cite{BeltranLeykin2011}. These constructions are similar to those given in 
Shub-Smale~\cite{Shub_Smale-V} and Armentano et al.~\cite{ABBCS} (this last, for the eigenvalue-eigenvector problem case).

The right-hand side in expression~\eqref{eq:bezsix} is known as 
the \emph{condition length} of the path $(f_t,\zeta_t)$. We will call \eqref{eq:bezsix} the {\em condition length estimate} of the number of steps.

Taking derivatives w.r.t. $t$ in the equality $f_t(\zeta_t)=0$ it is easily seen that
\begin{equation}\label{eq:zetadot}
   \dot \zeta_t=(Df_t(\zeta_t)|_{\zeta_{t}^\perp})^{-1} \dot f_t(\zeta_t),
\end{equation}
and with some work (see~\cite[Lemma 12, p.~231]{Blum_Cucker_Shub_Smale-1998} one can prove that
$$
  \|\dot\zeta\|_{\zeta_t}\leq \mu(f_t,\zeta_t) \|\dot f_t\|.
$$
It is known that $\mu(f,\zeta)\ge 1$, e,g., see \cite[Prop.~16.19]{Condition}.
So the estimate~\eqref{eq:bezsix} may be bounded from above by
\begin{equation}\label{eq:bezsixbound}
 K(f,g,\zeta)\leq C'\,D^{3/2}\, \int_0^1 \mu^2(f_t,\zeta_t)\,
 \|\dot f_t\|\,dt,
\end{equation}
where $C'=\sqrt{2} C$. Let us call this estimate the $\mu^2$-\emph{estimate}. 
\medskip
%\begin{remark}\label{rem:primerpaso}
%It is trivial to note from the proof of
%\eqref{eq:bezsix} that the following bound also holds:
%\begin{equation}\label{eq:bezsixb}
%   K=K(f,g,\zeta)\leq 1+ C'\,D^{3/2}\, \int_{t_1}^1 \mu^2(f_t,\zeta_t)\,
% \|\dot f_t\|\,dt,
%\end{equation}
%where $t_1$ is the size of the first homotopy step, that is $t_1=\Delta t_0$.
%\end{remark}

The condition length estimate is better than the $\mu^2$-estimate, but
algorithms achieving the smaller number of steps are more subtle and
the proofs of correctness more difficult.
 
Indeed in Beltr\'an-Pardo~\cite{Beltran_Pardo-2009-1} and
B\"urgisser-Cucker~\cite{Buergisser_Cucker-2010} the authors rely on
the $\mu^2$-estimate. At the times of these papers the algorithms
achieving the condition length bound where in development,
and~\cite{Buergisser_Cucker-2010} includes a construction which
achieves the $\mu^2$-estimate.

Yet, in a random situation, one might expect the improvement to be
similar to the improvement given by the average of $\|A(x)\|$, in all
possible directions, compared with $\|A\|$ (here, $A\colon\C^n\to\C^n$
denotes a linear operator), which according to
Armentano~\cite{Armentano} should give an improvement by a factor of
the square root of the domain dimension.  We have accomplished this
for the eigenvalue-eigenvector problem in Armentano et
al.~\cite{ABBCS}. Here we use an argument similar to that
of~\cite{ABBCS} to improve the estimate for the randomized algorithm
in Beltr\'an-Pardo~\cite{Beltran_Pardo-2011}.
 
The Beltr\'an-Pardo \texttt{randomized algorithm} works as follows
(see Beltr\'an-Pardo~\cite{Beltran_Pardo-2011}, and also
B\"urgisser-Cucker~\cite{Buergisser_Cucker-2010}): on input $f\in\hd$,
\begin{enumerate} 
\item 
Choose $f_0$ at random and then a zero $\zeta_0$ of $f_0$ at random. 
Beltr\'an and Pardo~\cite{Beltran_Pardo-2011} describe a 
general scheme to do so (roughly speaking, one first draws the ``linear'' 
part of $f_0$, computes $\zeta_0$ from it, and then draws the ``nonlinear'' 
part of $f_0$). An efficient implementation of this scheme, 
having running time $\CO(nDN)$, is fully described and analyzed  
in~\cite[Section~17.6]{Condition}.
\item 
Then connect $f_0/\|f_0\|$ to $f/\|f\|$ by an arc of a great
circle in the sphere, and invoke the continuation strategy above.
\end{enumerate}
The main result of~\cite{Beltran_Pardo-2011} is that the average number of
steps of this procedure is bounded by $O(D^{3/2}nN)$, and its total
average complexity is then $\CO(D^{3/2}nN^2)$ (since the
cost of an iteration of Newton's method, assuming all $d_i\geq2$, 
is $\CO(N)$, see~\cite[Proposition~16.32]{Condition}).

Our first main result is the following improvement of this last bound.

\begin{theorem}[{Randomized algorithm}]\label{thm:thm1}
The average number of steps of the randomized algorithm
with the condition length estimate is bounded by 
$$
  CD^{3/2}nN^{1/2} ,  
$$
where $C$ is a universal constant.
\end{theorem}

The constant $C$ can be taken as $\frac{\pi}{\sqrt{2}}C'$ with 
$C'$ not more than 400 even accounting for input and
round-off error, 
cf.\  Dedieu-Malajovich-Shub~\cite{Dedieu_Malajovich_Shub}.

\begin{remark}
Theorem~\ref{thm:thm1} is an improvement by a factor of $1/N^{1/2}$ of
the bound in~\cite{Beltran_Pardo-2011}, which results from using the
condition length estimate in place of the $\mu^2$-estimate.
\end{remark}

% 
% \noindent\textbf{Remark:}
% In Beltr\'an-Pardo~\cite{Beltran_Pardo-2011} and
% B\"urgisser-Cucker~\cite{Buergisser_Cucker-2010} the estimates on the
% complexity are $\mathcal O(D^{3/2} nN)$ in place of
% $\mathcal{O}(D^{3/2} nN^{1/2})$.  So
% the main theorem   The results of 
% Beltr\'an-Pardo~\cite{Beltran_Pardo-2011} gave a randomized solution
% to Smale's 17th problem, which asks if an approximation to one zero of
% a polynomial system $f$ in $\CH_{(d)}$ may be found in average
% polynomial time in the input size (see Smale~\cite{Smale-2000},
% Beltr\'an-Pardo~\cite{Beltran_Pardo-2011} and
% B\"urgisser-Cucker~\cite{Buergisser_Cucker-2010}).
  
Before proceeding with the proof of Theorem~\ref{thm:thm1}, 
we introduce some useful notation. We define the {\em solution variety} 
\[
   \V:=\{(f,\zeta)\in\hd\times\P(\C^{n+1})\mid f(\zeta)=0\},
\]
and consider the projections
\begin{equation}%\label{eq:projections}
%\begin{center}
\begin{tikzpicture}
\path (-0.25,0.3) node[right]{$\V$};
\draw[->] (0.15,0) -- (1.15,-1) node[above =3mm]{$\pi_2$};
\draw[->] (-0.15,0) -- (-1.15,-1) node[above =3mm]{$\pi_1$};
\path (-1.9,-1.35) node[right]{$\hd$};
\path (0.7,-1.35) node[right]{$\P(\C^{n+1})$.};
\end{tikzpicture}
\label{diagram:V}
% \caption{diagram:V}
%\end{center}
\end{equation}

The set of {\em ill-posed pairs} is the subset 
$$
  \Sigma':=\{(f,\zeta)\in\V\mid Df(\zeta)|_{\zeta^\perp} 
  \mbox{ is not invertible}\}= \{(f,\zeta)\in\V\mid\mu(f,\zeta) =\infty\}
$$
and its projection $\Sigma:=\pi_1(\Sigma')$ is the set of 
{\em ill-posed systems}. The number of 
iterations of the homotopy algorithm, $K(f,g,\zeta)$, is finite if and only if the lifting 
$\{(f_t,\zeta_t)\}_{t\in[0,1]}$ of the segment 
$\{f_t\}_{t\in[0,1]}$ does not cut $\Sigma'$.

\section{Proof of Theorem~\ref{thm:thm1}}

\subsection{Preliminaries}

Let us start this section with a few general facts we will use from
Gaussian measures.

Given a finite dimensional real vector space $V$ of dimension $m$,
with an inner product, we define two natural objects.
\begin{itemize}
\item The unit sphere $\S(V)$ with the induced Riemannian structure and 
volume form: the volume of $\S(V)$ is 
$\frac{2\pi^{m/2}}{\Gamma(\frac{m}{2})}$.

\item The Gaussian measure centered at $c\in V$, with variance $\frac{\sigma^2}{2}>0$, whose density is
\begin{equation}\label{def:gaussiandensity}
   \frac{1}{\sigma^m\pi^{m/2}}e^{-\|x-c\|^2/\sigma^2}.
\end{equation}

\end{itemize}

We will denote by ${N}_{V}(c,\sigma^2 \Id)$ the density given in
(\ref{def:gaussiandensity}). We will skip the notation of the
underlying space when it is understood. Furthermore, we will denote by
$\E_{x\in V}$ the average in the case $\sigma=1$ (that is, variance $1/2$).

%The standard gaussian, and the uniform  and these two objects are related as follows. 

\begin{lemma}\label{lem:polarcoord}
If $\varphi\colon V\to [0,+\infty]$ is measurable and homogeneous of degree $p>-m$, then
$$
 \E_{x\in V}(\varphi(x))=\frac{\Gamma(\frac{m+p}{2})}{\Gamma(\frac{m}{2})}
 \E_{u\in \S(V)}(\varphi(u)),
$$
where 
$$
 \E_{u\in \S(V)}(\varphi(u))=\frac{1}{\vol(\S(V))}\int_{\S(V)} \varphi(u)\,du.
$$
\end{lemma}

\proof
Integrating in polar coordinates we have
\begin{align*}
\E_{x\in V}(\varphi(x)) &
   =\frac{1}{\pi^{m/2}}\int_{x\in V}\varphi(x)\,e^{-\| x\|^2}\,dx\\
&=\frac{1}{\pi^{m/2}}\int_0^{+\infty}\rho^{m+p-1}e^{-\rho^2}\,d\rho
     \cdot\int_{u\in\S(V)}\varphi(u)\, du\\
&= \frac{\Gamma(\frac{m+p}{2})}{2\pi^{m/2}}\int_{u\in\S(V)}\varphi(u)\,du 
  =\frac{\Gamma(\frac{m+p}{2})}{\Gamma(\frac{m}{2})} 
  \E_{u\in\S(V)} (\varphi(u))
\end{align*} 
where we have used that $\int_0^{+\infty}\rho^{k}e^{-\rho^2}\,d
\rho=\frac12\Gamma(\frac{k+1}{2})$.
\eproof

The next results follows immediately from Fubini's theorem.

\begin{lemma}\label{lem:ortpro}
Let $E$ be a linear subspace of $V$, and let $\Pi:V\to E$ be the
orthogonal projection. Then, for any integrable function $\psi:E\rightarrow\R$ and for any	
$c\in V$, $\sigma>0$, we have
\begin{equation}\tag*{\qed}
   \E_{x\sim {N}_V(c,\sigma^2\Id)}(\psi(\Pi(x)))=\E_{y\sim {N}_E(\Pi(c),\sigma^2\Id)}(\psi (y)).
\end{equation}
\end{lemma}

When $V$ is a finite dimensional Hermitian vector space of complex
dimension~$m$, then the 
{\em complex Gaussian measure on $V$} with variance $\sigma^2$
is defined by the real Gaussian measure with variance $\sigma^2/2$ of the $2m$-dimensional 
real vector  space associated to $V$ with inner product the 
real part of the Hermitian product.

In this fashion, for any fixed $g\in\hd$ and $\sigma>0$, 
the Hermitian space
$(\hd,\pes\cdot\cdot)$ is equipped with the complex Gaussian measure
${N}(g,\sigma^2\Id)$. The expected value of a function
$\phi:\hd\to\R$ with respect to this measure is given by
\begin{equation}\label{def:GausMeas}
%\E_{f\in\hd}(\phi)=\frac{1}{\pi^N}\int_{f\in\hd}\phi(f)e^{-\|f\|^2}\,df.
\E_{f\sim {N}(g,\sigma^2\Id)}(\phi)=\frac{1}{\sigma^{2N}\pi^N}\int_{f\in\hd}\phi(f)e^{-\|f-g\|^2/\sigma^2}\,df.
\end{equation}

Fix any $\zeta\in\P(\C^{n+1})$. Following~\cite[Sect.~16.3]{Condition}, 
the space $\hd$ is orthogonally decomposed into the sum 
$C_\zeta\oplus\V_{\zeta}$, where 
\[ 
  C_{\zeta}=\bigg\{\diag\bigg(\czeta\bigg)a:\,
  a\in\C^n\bigg\}, 
\]
and 
\[
  \V_\zeta=\pi_2^{-1}(\zeta)=\{f\in\hd:\,f(\zeta)=0\}. 
\]
Note that $C_\zeta$ and $\V_\zeta$ are linear subspaces of $\hd$
of respective (complex) dimensions $n$ and $N-n$.  Note 
also that 
$$ 
   f_0=f-\diag\bigg(\czeta\bigg)f(\zeta) 
$$ 
is the orthogonal projection $\Pi_\zeta(f)$ of $f$ onto 
the fiber $\V_\zeta$.  This follows from
the reproducing kernel property of the Weyl Hermitian product on
$\CH_{d_i}$, namely,
\begin{equation}\label{eq:RKP}
\pes{g}{\pes{\cdot}{\zeta}^{d_i}}=g(\zeta),
\end{equation}
for all $g\in\mathcal{H}_{d_i}$ and $i=1,\ldots,n$.
In particular, the norm of $\pes{\cdot}{\zeta}^{d_i}\in\CH_{d_i}$ 
is equal to $\|\zeta\|^{d_i}$.

\subsection{Average condition numbers}

In this section we revisit the average value of the operator and Frobenius
condition numbers on $\hd$. The {\em Frobenius condition number} of $f$ at 
$\zeta$ is given by
\begin{equation}\label{eq:F-cn}
  \mu_F(f,\zeta) :=\|f\| \left\|\left(Df(\zeta)|_{\zeta^\perp} \right)^{-1}
   \diag(\|\zeta\|^{d_i-1}d_i^{1/2}) \right\|_F,
\end{equation}
that is, $\mu_F$ is defined as $\mu$ but using Frobenius instead 
of operator norm. Note that $\mu\leq\mu_F\leq\sqrt{n}\,\mu$.

Given $f\in\hd\setminus\Sigma$, the average of the condition numbers 
over the fiber is
$$
  \mu_{\mathrm{av}}^2(f):=\frac{1}{\mathcal D}\sum_{\zeta:\,f(\zeta)=0} 
  \mu^2(f,\zeta),\qquad \mufa^2(f):=
  \frac{1}{\mathcal D}\sum_{\zeta:\,f(\zeta)=0} \mu_{F}^2(f,\zeta)
$$
(or $\infty$ if $f\in\Sigma$). For simplicity, in what follows we write 
$\S:=\S(\CH_{(d)})$.  

Estimates on the probability distribution of the condition number
$\mu$ are known since~\cite{Bez2}. The exact expected value of $
\mu_{\mathrm{av}}^2(f)$ when $f$ is in the sphere $\S$ was found
in~\cite{Beltran_Pardo-2011} and the following estimate for the
expected value of $ \mu_{\mathrm{av}}^2(f)$ when $f$ is non-centered
Gaussian was proved in~\cite{Buergisser_Cucker-2010}: for all
$\hf\in\hd$ and all $\sigma>0$,
\begin{equation}\label{thm:MTR}
     \E_{f\sim N(\hf,\sigma^2\Id)} 
     \frac{\mu_{\av}^2(f)}{\|f\|^2}\ \leq\ \frac{e(n+1)}{2\sigma^2}.
\end{equation}

The following result slightly improves \eqref{thm:MTR}, even though it
is computed for~$\mu_F$.

\begin{theorem}[Average condition number]\label{theo:muFav} 
For every $\hf\in \hd$ and $\sigma>0$,
$$
 \E_{f\sim N(\hf,\sigma^2\Id)} 
     \frac{\mu_{F,\av}^2(f)}{\|f\|^2}\ \leq\  \frac{n}{\sigma^2},
$$
and equality holds in the centered case.
\end{theorem}

\begin{remark}\label{rmk:S}
The equality of Theorem~\ref{theo:muFav} implies from 
Lemma~\ref{lem:polarcoord} with $p=-2$ that
\[
   \E_{f\in\S} \mu_{F,\av}^2(f) =(N-1)n.
\]
\end{remark}

\begin{remark}
In the proof of Theorem~\ref{theo:muFav} we use the double-fibration
technology, a strategy based on the use of the %(modern) 
classical Coarea Formula, see for example
\cite[p.~241]{Blum_Cucker_Shub_Smale-1998}. In order to integrate some
real-valued function over $\hd$ whose value at some point $f$ is an
average over the fiber $\pi_1^{-1}(f)$, we lift it to $\V$ and then
pushforward to $\mathbb{P}(\C^{n+1})$ using the projections given in
(\ref{diagram:V}). The original expected value in $\hd$ is then writen
as an integral over $\mathbb{P}(\C^{n+1})$ which involves the quotient
of normal Jacobians of the projections $\pi_1$ and $\pi_2$. More
precisely,
\begin{equation}\label{eq:gendoublefibration}
  \int_{f\in\hd}\sum_{\zeta:\,f(\zeta)=0}\phi(f,\zeta)\,df
 =\int_{\zeta\in\mathbb{P}(\C^{n+1})}
   \int_{(f,\zeta)\in\pi_2^{-1}(\zeta)} \phi(f,\zeta) 
   \frac{\NJ_{\pi_1}}{\NJ_{\pi_2}}(f,\zeta)\,d\pi_2^{-1}(\zeta)\,d\zeta,
\end{equation}
where
\begin{equation*}
 \frac{\NJ_{\pi_1}}{\NJ_{\pi_2}}(f,\zeta)
 = |\det(Df(\zeta)|_{\zeta^\perp})|^2
\end{equation*}
(see~\cite[Section~13.2]{Blum_Cucker_Shub_Smale-1998},%
~\cite[Section~17.3]{Condition}, or~\cite[Theorem 6.2]{ABBCS} 
for further details and other examples of use).

We point out that the proof of Theorem~\ref{theo:muFav} can also be achieved using the (slightly) different method of~\cite{Beltran_Pardo-2011} and \cite[Chapter 18]{Condition} based on the mapping taking $(f,\zeta)$ to $(Df(\zeta),\zeta)$ whose Jacobian is known to be constant (see~\cite[Main Lemma]{Beltran_Pardo-2011}).
\end{remark}

\proofof{Theorem~\ref{theo:muFav}}
By the definition of non-centered Gaussian, and the double-fibration 
formula~(\ref{eq:gendoublefibration}), we have
 \begin{align}\label{eq:start-formula}
& \E_{f\sim N(\hf,\sigma^2\Id)} 
     \frac{\mu_{F,\av}^2(f)}{\|f\|^2} =\frac{1}{\mathcal{D}} \int_{f\in\hd} 
  \bigg(\sum_{\zeta:\,f(\zeta)=0}\frac{\mu_F^2(f,\zeta)}{\|f\|^2}\bigg)\,
  \frac{e^\frac{-\|f-\hf\|^2}{\sigma^2}}{\sigma^{2N}\pi^N}\,df \\ \notag
 &=\frac{1}{\mathcal{D}} \frac{1}{\sigma^{2n}\pi^n}\int_{\zeta\in\P(\C^{n+1})}
e^\frac{-\|\hf-\Pi_{\zeta}(\hf)\|^2}{\sigma^2}
 \int_{f\in \V_\zeta}\frac{\mu_F^2(f,\zeta)}{\|f\|^2}\,
   \big|\det(Df(\zeta)|_{\zeta^\perp}\big|^2\,
\frac{e^\frac{-\|f-\Pi_{\zeta}(\hf)\|^2}{\sigma^2}}{(\sigma^2\pi)^{N-n}}\,
df\,d\zeta,
\end{align}
where we have used that
$\|f-\hf\|^2=\|f-\Pi_{\zeta}(\hf)\|^2+\|\hf-\Pi_{\zeta}(\hf)\|^2$
for every $f\in\V_\zeta$. 

We simplify now the integral $I_\zeta$ over the fiber $\V_\zeta$. 
Let $U_\zeta$ be a unitary transformation of $\C^{n+1}$ such that
$U_\zeta(\zeta/\|\zeta\|)=e_0$. Then, by the invariance  under unitary
substitution of each term under the integral sign, 
we have by the change of variable formula that 
\begin{align*}
I_\zeta\, :=\,& \int_{f\in \V_\zeta}\frac{\mu_F^2(f,\zeta)}{\|f\|^2}\,
   \big|\det(Df(\zeta)|_{\zeta^\perp}\big|^2\,
\frac{e^\frac{-\|f-\Pi_{\zeta}(\hf)\|^2}{\sigma^2}}{\sigma^{2(N-n)}\pi^{N-n}}\, df \\
=\,& 
\int_{h\in \V_{e_0}}\frac{\mu_F^2(h,e_0)}{\|h\|^2}\,
   \big|\det(Dh(e_0)|_{e_0^\perp}\big|^2\,
\frac{e^\frac{-\|h-\Pi_{e_0}(\widehat{h})\|^2}{\sigma^2}}{\sigma^{2(N-n)}\pi^{N-n}}\, dh,
\end{align*}
where $\widehat{h}_\zeta:= \hf\circ U_{\zeta}^{-1}$. 
We project now  $h\in\V_{e_0}$ orthogonally onto the vector space
\[
L_{e_0} :=\{g\in\hd: g(e_0)=0,D^kg(e_0)=0\text{ for }k\geq2\} ,
\]
obtaining $g\in L_{e_0}$. Since 
$Dh(e_0)|_{e_0^\perp}$ coincides with $Dg(e_0)|_{e_0^\perp}$ 
(see for example~\cite[Prop.~16.16]{Condition}), we conclude 
by Fubini that
\begin{align*}
I_\zeta \,=\,& \int_{h\in \V_{e_0}}\frac{\mu_F^2(h,e_0)}{\|h\|^2}\,
   \big|\det(Dh(e_0)|_{e_0^\perp}\big|^2\,
\frac{e^\frac{-\|h-\Pi_{e_0}(\widehat{h}_\zeta)\|^2}{\sigma^2}}{\sigma^{2(N-n)}\pi^{N-n}}\,dg\\
=\,& 
 \int_{g\in L_{e_0}}\frac{\mu_F^2(g,e_0)}{\|g\|^2}\,
   \big|\det(Dg(e_0)|_{e_0^\perp}\big|^2\,
\frac{e^\frac{-\|g-\Pi_{L_{e_0}}(\widehat{g}_\zeta)\|^2}{\sigma^2}}{\sigma^{2n^2}\pi^{n^2}}\,
dh,
\end{align*}
where $\widehat{g}_\zeta := \Pi_{L_0}(\widehat{h}_\zeta)$. 
%and $\Pi_{L_{e_0}}$ is the orthogonal projection onto $L_{e_0}$. 
By the change of variable given by
$$
 L_{e_0}\to \C^{n\times n},\  g\mapsto A 
 :=\diag(d_i^{-1/2})Dg(e_0)|_{e_0^\perp} ,
$$
we have $\frac{\mufa^2(g)}{\|g\|^2}=\|A^{-1}\|_F^2$ 
and denoting by $\widehat{A}_\zeta$ the image of $\widehat{g}_\zeta$, 
we obtain that 
\begin{equation*}\label{eq:def-I}
 I_\zeta \ = \E_{A\in{N}(\widehat{A}_\zeta,\sigma^2\Id_n)} \|A^{-1}\|_F^2\,
 |\det(A)|^2 .
\end{equation*}
We thus conclude from \eqref{eq:start-formula} that 
\begin{equation}\label{eq:uno}
 \E_{f\sim N(\hf,\sigma^2\Id)} \frac{\mufa^2(f)}{\|f\|^2} 
  =\  \frac{1}{\mathcal{D}} \frac{1}{\sigma^{2n}\pi^n}\int_{\zeta\in\P(\C^{n+1})} 
 e^\frac{-\|\hf-\Pi_{\zeta}(\hf)\|^2}{\sigma^2} \, I_\zeta \, d\zeta. 
\end{equation}
If we replace $\mu_{F,\av}(f)^2/\|f\|^2$ by the constant function $1$ on $\hd$, 
the same argument leading to~\eqref{eq:uno} now leads to 
\begin{equation}\label{eq:eins}
 1 %= \E_{f\sim N(\hf,\sigma^2\Id)}\big(\ 1 \big) 
 \ =\  \frac{1}{\mathcal{D}} \frac{1}{\sigma^{2n}\pi^n}\int_{\zeta\in\P(\C^{n+1})} 
 e^\frac{-\|\hf-\Pi_{\zeta}(\hf)\|^2}{\sigma^2} \, 
   \E_{A\in{N}(\widehat{A}_\zeta,\sigma^2\Id_n)} |\det(A)|^2 \, d\zeta. 
\end{equation}

From Proposition 7.1 of Armentano et al.~\cite{ABBCS}, we can bound 
\begin{equation}\label{eq:diegos-bd}
   I_\zeta \ = \E_{A\in{N}(\widehat{A}_\zeta,\sigma^2\Id_n)} \|A^{-1}\|_F^2\, |\det(A)|^2  
  \ \leq\  \frac{n}{\sigma^2}  \E_{A\in{N}(\widehat{A}_\zeta,\sigma^2\Id_n)} |\det(A)|^2 
\end{equation}
(with equality if $\widehat{A}_\zeta=0$).
By combining \eqref{eq:uno}, \eqref{eq:diegos-bd}, and \eqref{eq:eins}
we obtain
$$
\E_{f\sim N(\hf,\sigma^2\Id)} \frac{\mufa^2(f)}{\|f\|^2}\ \leq\ \frac{n}{\sigma^2} .
$$
Moreover, equality holds if $\hf=0$ and hence $\widehat{A}_\zeta=0$
for all $\zeta$.
\eproof

\subsection{Complexity of the randomized algorithm}

The goal of this section is to prove Theorem~\ref{thm:thm1}. 
To do so, we begin with some preliminaries.

For $f\in\S$ we denote by $T_f\S$ the tangent space at $f$ of $\S$. 
This space is equipped with the real part of the Hermitian structure of $\hd$,
and coincides with the (real) orthogonal complement of $f\in\CH_{(d)}$. 

We consider the map 
$\phi\colon\S\times \hd \to [0,\infty]$ 
defined for $f\not\in\Sigma$ by 
\[
  \phi(f,\dot f):=\frac{1}{\mathcal{D}}\sum_{\zeta:\,f(\zeta)=0}
  \mu(f,\zeta)\big\|(\dot f,\dot\zeta\big)\| ,
\]
where $\dot \zeta = (Df(\zeta)|_{\zeta^\perp})^{-1} \dot{f} (\zeta)$, 
and by  $\phi(f,\dot{f}) := \infty$ if $f\in\Sigma$.
Note that $\phi$ satisfies 
$\phi(f,\lambda \dot{f}) = \lambda\phi(f,\dot{f})$ 
for $\lambda\ge 0$. 

Suppose that $f_0,f\in\S$ are such that $f_0\ne\pm f$ and 
denote by $\mathcal{L}_{f_0,f}$ the shorter great circle segment  
with endpoints $f_0$ and $f$. 
Moreover, let 
$\a = d_\S(f_0,f)$ denote the angle between $f_0$ and $f$. 
If  $[0,1]\to\S,\, t\mapsto f_t$ is the constant speed 
parametrization of $\mathcal{L}_{f_0,f}$ with endpoints $f_0$ and $f_1=f$, 
then $\|\dot{f}_t\|=\a$. We may also parametrize $\mathcal{L}_{f_0,f}$ 
by the arc-length $s=\a t$, setting $F_s:=f_{\a^{-1}s}$, in which case 
$\dot{F}_s = \a^{-1}\dot{f}_t$ is the unit tangent vector (in the 
direction of the parametrization)
to $\mathcal{L}_{f_0,f}$ at $F_s$. Moreover, 
$$
 \int_0^1\phi(f_t,\dot  f_t)\,dt 
  =  \int_0^{\a}\phi(F_s,\dot{F}_s)\,ds .
$$
Consider the compact submanifold $\mathcal{S}$ of $\S\times\S$ given by 
$$ 
   \CS = \{(f,\dot f)\in\S\times\S:\dot f\in T_f\S\},
$$ 
where $T_f\S$ denotes the real tangent space of $\S$ at $f$. 
We endow $\S$ with the Riemannian metric induced from the real part 
of the Hermitian product of $\hd$, and therefore $\mathcal{S}$ inherits 
the product Riemannian structure.

The following lemma has been proven in Armentano et al.~\cite{ABBCS}.

\begin{lemma}\label{lemn:sphere}
We have 
$$
 I_\phi := \E_{f_0,f\in\S}\left(\int_0^{1} 
     \phi(f_t,\dot{f}_t)\,dt\right) 
 =  \frac{\pi}{2}\, \E_{(f,\dot f)\in\CS}  
      \big(\phi(f,\dot{f}) \big) ,
$$
where the expectation on the right hand-side refers to 
the uniform distribution on~$\CS$. 
\end{lemma}

We proceed with a further auxiliary result. 
For $f\in\S$ we consider the unit sphere
$\CS_f:=\{\dot{f} \in T_f\S :(f,\dot f)\in\CS\}$  in $T_f\S$.

\begin{lemma}\label{lem:l4}
Fix $f\in\S$ and $\zeta\in\P(\C^{n+1})$ 
with $f(\zeta)=0$. 
For $\dot{f} \in \CS_f$  
let $\dot \zeta$ be the function of $(f,\dot{f})$ 
and $\zeta$ given by 
$\dot \zeta = (Df(\zeta)|_{\zeta^\perp})^{-1} \dot{f} (\zeta)$. 
Then we have
$$
  \E_{\dot f\in \CS_f}(\|\dot\zeta\|^2)
  = \frac{1}{N-\frac12}\, \big\|(Df(\zeta)|_{\zeta^\perp})^{-1}\big\|_F^2 ,
$$
where the expectation is with respect to the uniform probability distribution 
on $\CS_f$. 
\end{lemma}

\proof
Since the map $T_f\S\to\R,\, \dot{f} \mapsto \|\dot{\zeta}(\dot{f})\|^2$ 
is quadratic, we get from Lemma~\ref{lem:polarcoord} 
(recall that $\dim T_f\S = 2N -1$)
$$
  \E_{\dot f\in T_f\S}(\|\dot\zeta(\dot f)\|^2) 
 = {\Big(N-\frac12\Big)}\, \E_{\dot f\in\mathcal{S}_f} 
 \big(\big\|\dot\zeta(\dot f)\big\|^2\big).
$$
Note that the mapping $\CH_{(d)}\to C_\zeta$ given by $\dot f\mapsto
\Pi_{C_\zeta}\dot f$ is an orthogonal projection, and furthermore
$C_\zeta\to \C^n$ given by $\dot f\mapsto \dot f(\zeta)$ is a linear
isometry. Then from Lemma~\ref{lem:ortpro}, and the change of 
variables formula we obtain 
$$
  \E_{\dot f\in T_f\S} \big(\big\|\dot \zeta(\dot f)\big\|^2\big)
 =\E_{\dot w\in{\C^n}} \big\|(Df(\zeta)|_{\zeta^\perp})^{-1}\dot w\big\|^2
 = \big\|(Df(\zeta)|_{\zeta^\perp})^{-1}\big\|_F^2,
$$
where the last equality is straightforward. 
\eproof

\noindent
\proofof{Theorem~\ref{thm:thm1}}
From~\eqref{eq:bezsix}, using the notation from there, 
we know that the number of Newton steps 
of the homotopy with starting pair $(f_0,\zeta_0)$ and target system~$f$ 
is bounded as 
$$
 K(f,f_0,\zeta_0) \ \leq\ 
  C D^{3/2} \, \int_0^1 \mu(f_t,\zeta_t)\, \| (\dot{f}_t,\dot{\zeta}_t)\|\, dt .
$$
Hence we get for $f,f_0\in\S$, 
\begin{align*}
  \frac{1}{\mathcal D} \sum_{\zeta_0:\,f_0(\zeta_0)=0} K(f,f_0,\zeta_0) 
  &\ \leq
   C D^{3/2} \, \int_0^1 \frac{1}{\mathcal{D}} \sum_{\zeta_0:\,f_0(\zeta_0)=0}
    \mu(f_t,\zeta_t)\, \| (\dot{f}_t,\dot{\zeta}_t)\|\, dt  \\
   &\ = 
    C D^{3/2} \, \int_0^1 \phi(f_t,\dot{f}_t)\, dt . 
\end{align*}
Therefore, by Lemma~\ref{lemn:sphere}, 
\begin{equation}\label{eq:BPconleness}
  \E_{f,f_0\in\S} 
  \Big(\frac{1}{\mathcal{D}}\sum_{\zeta_0:\,f_0(\zeta_0)=0}
    K(f,f_0,\zeta_0)\Big)
   \ \leq\ C\,D^{3/2} \, \frac{\pi}{2}\, \E_{(f,\dot f)\in\mathcal{S}}  
      \big(\phi(f,\dot{f}) \big) .
\end{equation}
It is easy to check that the projection 
$\CS\to\S,\, (f,\dot{f})\mapsto f$, 
has the Normal Jacobian $1/\sqrt{2}$.
From the coarea formula, we therefore obtain 
\begin{align*}
  \E_{(f,\dot f)\in\mathcal{S}} \big( \phi(f,\dot f) \big)
  &= \sqrt{2}\E_{f\in\S}\ \E_{\dot{f}\in \CS_f} \big( \phi(f,\dot f)  \big) \\
  &=\sqrt2\, \E_{f\in\S}\Big(\frac{1}{\mathcal{D}}
  \sum_{\zeta:\,f(\zeta)=0}\mu(f,\zeta)\E_{\dot f\in \mathcal{S}_f}
  \big(\big\|(\dot f,\dot\zeta)\big\|\big)\Big).
\end{align*}
In order to estimate this last quantity, note first that from the 
Cauchy-Schwartz inequality
\begin{eqnarray*}
  \E_{\dot f\in \mathcal{S}_f}
 \big(( 1+ \|\dot{\zeta}\|^2)^{\frac12} \big)
  &\leq &\Big(1+\E_{\dot f\in \mathcal{S}_f}
 (\|\dot\zeta\|^2)\Big)^{1/2}\\
&\leq& \left(1+\frac{1}{N-\frac12}\, 
   \big\|(Df(\zeta)|_{\zeta^\perp})^{-1}\big\|_F^2\right)^{1/2}
\end{eqnarray*}
the last by Lemma~\ref{lem:l4}. 
Now we use 
$\|(Df(\zeta)|_{\zeta^\perp})^{-1}\|_F\leq \mu_F(f,\zeta)$ 
and $\mu(f,\zeta)\leq \mu_F(f,\zeta)$ to deduce  
\begin{align*}
\frac{1}{\sqrt{2}} \E_{(f,\dot f)\in\mathcal{S}} \big( \phi(f,\dot f) \big) 
  & \leq\ 
  \E_{f\in\S}\bigg(\frac{1}{\mathcal{D}}
  \sum_{\zeta:\,f(\zeta)=0}\mu_F(f,\zeta)
  \Big(1+\frac{ \mu_F^2(f,\zeta)}{N-\frac12}\Big)^\frac{1}{2}\bigg)\\
 &\ \leq\ 
  \E_{f\in\S}\bigg(\frac{1}{\mathcal{D}}
  \sum_{\zeta:\,f(\zeta)=0}\bigg( \frac{(N-\frac12)^{\frac12}}{2}
 +\frac{ \mu_F^2(f,\zeta)}{(N-\frac12)^{\frac12}}\bigg)\bigg)\\
 &\ =\ 
  \frac{(N-\frac12)^{\frac12}}{2}+ \E_{f\in\S}\bigg(
  \frac{ \mu_{F,\av}^2(f)}{(N-\frac12)^{\frac12}}\bigg)
\end{align*}
the second inequality since for all $x\ge 0$ and $a>0$ 
we have 
$$
 x^{1/2}(1+a^2x)^{1/2}\leq \frac{1}{2a}+ax.
$$
A call to Remark~\ref{rmk:S} finally yields 
\[
  \frac{1}{\sqrt{2}} \E_{(f,\dot f)\in\mathcal{S}} \big( \phi(f,\dot f) \big)
 \ \leq\ 
  \frac{(N-\frac12)^{\frac{1}{2}}}{2} +\frac{(N-1)n}{(N-\frac12)^{\frac12}} 
 \ \leq\  \sqrt{N}\bigg(\frac12+n\bigg) .
\]
Replacing this bound in~\eqref{eq:BPconleness} finishes 
the proof.\eproof

\section{A Deterministic Algorithm}

A deterministic solution for Smale's 17th problem is yet to be found. 
The state of the art for this theme is given 
in~\cite{Buergisser_Cucker-2010} where the following result is 
proven.

\begin{theorem}\label{th:near17}
There is a deterministic real-number algorithm that on 
input~$f\in\hd$ computes an approximate zero of~$f$ in 
average time $N^{\CO(\log\log N)}$.
Moreover, if we restrict data to polynomials satisfying
$$
    D\le n^{\frac1{1+\ep}} \quad\mbox{ or }\quad D\ge n^{1+\ep},
$$
for some fixed $\ep>0$,
then the average time of the algorithm is polynomial
in the input size~$N$.
\end{theorem}

The algorithm exhibited in~\cite{Buergisser_Cucker-2010} 
uses two algorithmic strategies according to whether $D\leq n$ or 
$D>n$. In the first case, it applies a homotopy method 
and in the second an adaptation of a method coming from 
symbolic computation. 

The goal of this section is to show that a more unified approach, 
where homotopy methods are used in both cases, yields a 
proof of Theorem~\ref{th:near17} as well. Besides a gain in 
expositional simplicity, this approach can claim for it the 
well-established numerical stability of homotopy methods. 

In all what follows we assume the simpler homotopy 
algorithm
in~\cite{Buergisser_Cucker-2010} 
(as opposed to those 
in~\cite{Beltran-2011,Dedieu_Malajovich_Shub}). Its choice of 
step length at the $k$th iteration is proportional 
to $\mu^{-2}(f_{t_k},x_{t_k})$ (which, in turn, is 
proportional to $\mu^{-2}(f_{t_k},\zeta_{t_k})$). 
For this algorithm, we have the 
$\mu^2$-estimate~\eqref{eq:bezsixbound} but not
the finer estimate~\eqref{eq:bezsix}.

To understand the technical requirements of the analysis of a
deterministic algorithm, let us summarize an analysis (simpler than
the one in the preceding section because of the assumption above) for
the randomized algorithm.  Recall, the latter draws an initial pair
$(g,\zeta)$ from a distribution which amounts to first draw $g$ from
the distribution on $\S$ and then draw $\zeta$ uniformly among the
$\CD$ zeros $\{\zeta^{(1)},\ldots,\zeta^{(\CD)}\}$ of $g$. The
$\mu^2$-estimate~\eqref{eq:bezsixbound} provides an upper bound for
the number of steps needed to continue $\zeta$ to a zero of $f$
following the great circle from $g$ to $f$ (assuming $\|f\|=\|g\|=1$
and $f\neq \pm g$). Now~\eqref{eq:bezsixbound} does not change if we
reparametrize $\{f_t\}_{t\in [0,1]}$ by arc--length, so we can also write 
it as
\[
   K(f,g,\zeta)\leq C'\,D^{3/2}\, \int_0^{d_\S(g,f)} \mu^2(f_s,\zeta_s)\,ds,  
\]
where $d_\S(g,f)$ is the spherical distance from $g$ to $f$.
Thus, the average number of homotopy iterations 
satisfies
\begin{align}\notag
 \E_{f\in\S}\E_{g\in\S}\frac{1}{\CD}\sum_{i=1}^{\CD} K(f,g,\zeta^{(i)})
 &\leq\ C'\,D^{3/2}\, \E_{f\in\S}\E_{g\in\S} \frac{1}{\CD}\sum_{i=1}^{\CD}
  \int_0^{d_\S(g,f)} \mu^2(f_s,\zeta_s)\,ds \\ \label{eq:new}
 & \leq\  C'\,D^{3/2}\, \E_{f\in\S}\E_{g\in\S}  \int_0^{d_\S(g,f)} 
  \mu_{F,\av}^2(f_s)\,ds .
\end{align}
Let $P_s$ denote the set of pairs $(f,g) \in\S^2$ such that $d_\S(g,f) \ge s$. 
Rewriting the above integral using Fubini, we get 
$$
 \E_{f\in\S}\E_{g\in\S}  \int_0^{d_\S(g,f)} \mu_{F,\av}^2(f_s)\,ds  
 = \int_0^{\pi} \int_{P_s}  \mu_{F,\av}^2(f_s)\, dfdg\, ds  
= \frac{\pi}{2} \E_{h\in\S} \mu_{F,\av}^2(h),
$$
the second equality holding since 
for a fixed $s\in[0,\pi]$ and uniformly distributed $(f,g)\in P_s$,  
one can show that the system $f_s$ is uniformly distributed on $\S$. 
Summarizing, we get 
\[
  \E_{f\in\S}\E_{g\in\S}\frac{1}{\CD}\sum_{i=1}^{\CD} K(f,g,\zeta^{(i)})
 \ \leq\ C'\,D^{3/2}\, \frac{\pi}{2} \E_{h\in\S}\mufa^2(h)
 \underset{\text{Rmk.~\ref{rmk:S}}}{=}
 C'\,D^{3/2}\, \frac{\pi}{2} (N-1)n . 
\]
This constitutes an elegant derivation of 
the previous $\CO(nD^{3/2}N)$ bound (but not of the sharper bound of our Theorem~\ref{thm:thm1}).
\medskip

\proofof{Theorem~\ref{th:near17}}
If the initial pair $(g,\zeta)$ is not going to be random we face 
two difficulties. Firstly ---as $g$ is not random--- 
the intermediate systems $f_t$ are not going 
to be uniformly distributed on $\S$. Secondly
---as $\zeta$ is not random---
we will need a bound 
on a given $\mu^2(f_t,\zeta_t)$ rather than one on the mean  
of these quantities (over the $\CD$ possible zeros of $f_t$), as 
provided by Theorem~\ref{theo:muFav}.  

Consider a fixed initial pair $(g,\zeta)$ with $g\in\S$ and let 
$s_1$ be the step length of the first step of the algorithm (see
for example the definition of Algorithm {\rm ALH}
in~\cite{Buergisser_Cucker-2010}), which satisfies
\begin{equation}\label{eq:s1}
 s_1\geq \frac{c}{D^{3/2}\mu^2(g,\zeta)}\quad\text{ ($c$ a constant).}
\end{equation}
Note that this bound on the length $s_1$ of the first homotopy step
depends on the condition $\mu(g,\zeta)$ only and is thus independent
of the condition at the other zeros of $g$.  Any of the mentioned
versions of the continuation algorithm, not only that
of~\cite{Buergisser_Cucker-2010}, satisfies \eqref{eq:s1}. 

Consider also the (small) portion of great circle contained in $\S$ 
with endpoints $g$ and $f/\|f\|$, which we parametrize by arc length and
call $h_s$ (that is, $h_0=g$ and $h_\alpha=f/\|f\|$ where
$\alpha=d_\S(g,f/\|f\|)$) defined for $s\in[0,\alpha]$. Thus, 
after the first step of the homotopy, the current pair is 
$(h_{s_1},x_1)$ and we denote by $\zeta'$ the zero of $h_{s_1}$ 
associated to $x_1$.
For a time to come, we will focus on bounding the quantity 
$$
 H:=H(g,\zeta) 
 :=\E_{f\in\hd}\frac{1}{\CD}\sum_{i=1}^{\CD}
 K\big(f/\|f\|,h_{s_1},\zeta^{(i)}\big),
$$
where the sum is over all the zeros $\zeta^{(i)}$ of $h_{s_1}$. 
This is the average of the number of homotopy steps over {\em both}
the system $f$ and the $\CD$ zeros of $h_{s_1}$. 
We will be interested in this average even though we will not 
consider algorithms following a path randomly chosen: 
the homotopy starts at the pair $(g,\zeta)$, moves to $(h_{s_1},x_1)$ 
and proceeds following this path. 

From~\eqref{eq:bezsixbound} applied to $(h_{s_1},\zeta')$,
\begin{equation}\label{eq:K}
 K\big(f/\|f\|,h_{s_1},\zeta^{(i)}\big) 
\ \leq\ C'\,D^{3/2}\, 
 \int_{s_1}^{\alpha} \mu^2(h_s,\zeta^{(i)}_s)\, \|\dot h_s\|\,ds,
\end{equation}
Reparametrizing $\{h_s:s_1\leq s\leq\alpha\}$ by
$\{f_t/\|f_t\|:t_1\leq t\leq 1\}$ where $f_t=(1-t)g+tf$ and $t_1$ is
such that $f_{t_1}/\|f_{t_1}\|=h_{s_1}$ (see Lemma~\ref{lem:t1} below)
does not change the value of the path integral in
\eqref{eq:K}. Moreover, a simple computation shows that
\[
 \left\|\frac{d}{dt}\left(\frac{f_t}{\|f_t\|}\right)\right\|
 \ \leq\ \frac{\|f\|\|g\|}{\|f_t\|^2}=\frac{\|f\|}{\|f_t\|^2},
\]
so we have
\begin{equation}\label{eq:bezsix2}
  K\big(f/\|f\|,h_{s_1},\zeta^{(i)}\big)
 \ \leq\ C'\,D^{3/2}\, \|f\|\,\int_{t_1}^1 \frac{\mu^2(f_t,\zeta^{(i)}_t)}
 {\|f_t\|^2}\,dt.
\end{equation}
Because of scale invariance, the quantity $H$ satisfies  
$$
 H = \E_{f\in\hd^{\sqrt{2N}}}\frac{1}{\CD}\sum_{i=1}^{\CD}
 K(f,h_{s_1},\zeta^{(i)}) ,
$$
where the second expectation is taken over a truncated Gaussian 
(that only draws systems $f$ with $\|f\|\leq\sqrt{2N}$) with density 
function given by 
$$
  \rho(f):=\begin{cases}\frac1{p}\varphi(f) & \mbox{if $\|f\|\leq\sqrt{2N}$}\\
               0 & \mbox{otherwise.}  
               \end{cases}
$$
Here $\varphi$ is the density function of the standard Gaussian on 
$\hd$ and $p:=\Prob\{\|f\|\leq \sqrt{2N}\}$. 
Then, using~\eqref{eq:bezsix2}, 
\[
 H\ \leq\ 
   \sqrt{2N}C'\,D^{3/2}\,
  \E_{f\in\hd^{\sqrt{2N}}}\frac{1}{\CD}\sum_{i=1}^{\CD}
  \int_{t_1}^1\frac{\mu^2(f_t,\zeta_t^{(i)})}{\|\dot f_t\|^2}\,dt.
\]
From Lemma~\ref{lem:t1} we have
$$
  t_1 \geq \frac{c'}{D^{3/2}\sqrt{N}\mu^2(g,\zeta)}
$$  
for a constant $c'$ (different from, but close to, $c$). We thus 
have proved that there are constants $C'',c'$ 
such that
\begin{align*}
 H\ \leq\ &
   C''\sqrt{N}\,D^{3/2}\,
  \E_{f\in\hd^{\sqrt{2N}}}\frac{1}{\CD}\sum_{i=1}^{\CD}
  \int_{\frac{c'}{D^{3/2}\sqrt{N}\mu^2(g,\zeta)}}^1
 \frac{\mu^2(f_t,\zeta_t)}{\|\dot f_t\|^2}\,dt\\
  =\ &C''\sqrt{N}\,D^{3/2}\,
  \E_{f\in\hd^{\sqrt{2N}}}
  \int_{\frac{c'}{D^{3/2}  \sqrt{N}\mu^2(g,\zeta)}}^1
  \frac{\mu_{\av}^2(f_t)}{\|\dot f_t\|^2}\,dt.
\end{align*}
The $2N$ in the cut-off for the truncated Gaussian is the 
expectation of $\|f\|^2$ for a standard Gaussian $f\in\hd$. 
Since this expectation is not smaller than the median of 
$\|f\|^2$ (see~\cite[Cor.~6]{Choi}) we have 
$\frac{1}p\leq 2$. Using this inequality and the fact that 
the random variable we are taking the expectation of is nonnegative, 
we deduce that 
\begin{align*}
  H\ \leq\ & 
   2C''\sqrt{N}\,D^{3/2}\,
  \E_{f\in\hd}
  \int_{\frac{c'}{D^{3/2} \sqrt{N}\mu^2(g,\zeta)}}^1
  \frac{\mu_{\av}^2(f_t)}{\|\dot f_t\|^2}\,dt\\
  \leq\ &2C''\sqrt{N}\,D^{3/2}\,
  \int_{\frac{c'}{D^{3/2} \sqrt{N}\mu(g,\zeta)^2}}^1\E_{f\in\hd}
  \frac{\mufa^2(f_t)}{\|\dot f_t\|^2}\,dt.
\end{align*}
We next bound the expectation in the right-hand side 
using Theorem~\ref{theo:muFav} and the fact that 
$f_t\sim N((1-t)g,t^2\Id)$ and obtain 
\begin{eqnarray}\label{eq:homot3}
  H&\leq& 2nC''\sqrt{N}\,D^{3/2}\,
  \int_{\frac{c'}{D^{3/2} \sqrt{N}\mu^2(g,\zeta)}}^1
  \frac{1}{t^2}\,dt\nonumber\\
     &\leq& C'''\,D^3nN\mu^2(g,\zeta),
\end{eqnarray}
with $C'''$ yet another constant.

Having reached thus far, the major obstacle we face 
is that the quantity $H$, for which we derived the 
bound~\eqref{eq:homot3}, is an average over all 
initial zeros of $h_{s_1}$ (as well as over $f$). For this 
obstacle, none of the two solutions below is fully 
satisfactory, but each of them is so for a broad choice of 
pairs $(n,D)$. 
\medskip

\noindent
{\bf Case 1: $D>n$.\quad}
Consider any $g\in\S$, $\zeta$ a well-posed zero of $g$, 
and let $\zeta^{(1)},\ldots,\zeta^{(\CD)}$ 
be the zeros of $h_{s_1}$.  
Note that when $f$ is Gaussian, these are $\CD$ different zeros 
almost surely. Clearly,
\begin{eqnarray*}
  \E_{f\in\hd}  K(f,g,\zeta) &\leq& 
  1+ \E_{f\in\hd}\sum_{i=1}^{\CD}
  K(f,h_{s_1},\zeta^{(i)}) \\ 
 &=& 1+\CD\, H =\CO(\CD D^3Nn\mu^2(g,\zeta))
\end{eqnarray*}
the last by~\eqref{eq:homot3}. We now take as initial pair 
$(g,\zeta)$ the pair $(\og, e_0)$ where 
$\og=(\og_1,\ldots,\og_n)$ is given by 
$$
    \og_i=(d_1^{-1}+\cdots+d_n^{-1})^{-1/2} X_0^{d_i-1}X_i,
    \qquad\mbox{for $i=1,\ldots,n$}
$$
(the scaling factor guaranteeing that $\|\og\|=1$) 
and $e_0=(1,0,\ldots,0)\in\C^{n+1}$. It is known that 
$\mu(\og, e_0)=1$ 
(see~\cite[Rem.~16.18]{Condition})
(and that all 
other zeros of $\og$ are ill-posed, but this is not 
relevant for our argument). Replacing this equality 
in the bound above we obtain
\begin{equation}\label{eq: complexityboundC1}
\E_{f\in\hd}  K(f,\og,e_0) =\CO(\CD D^3Nn) ,
\end{equation}
which implies an average cost of 
$\CO(\CD D^3N^2n)$ since the number of operations 
at each iteration of the homotopy algorithm is $\CO(N)$
(see~\cite[Proposition~16.32]{Condition}). 

For any $\ep>0$ this quantity is polynomially bounded 
in $N$ provided $D\geq n^{1+\ep}$ and is bounded 
as $N^{\CO(\log\log N)}$ when $D$ is in the range 
$[n,n^{1+\ep}]$ (\cite[Lemma~11.1]{Buergisser_Cucker-2010}).  
\medskip

\noindent
{\bf Case 2: $D\leq n$.\quad}
The occurrence of $\CD$ makes the bound 
in~\eqref{eq: complexityboundC1} too large 
when $D$ is small. In this case, we consider 
the initial pair $(\oU,\bz_1)$ where $\oU\in\hd$ is given by  
\begin{equation*}\label{eq:def_g}
  \oU_1 = \frac1{\sqrt{2n}} ( X_0^{d_1} - X_1^{d_1}), \ldots,\, 
   \oU_n = \frac1{\sqrt{2n}} ( X_0^{d_n} - X_n^{d_n} ) ,
\end{equation*}
(the scaling factor guaranteeing that $\|\oU\| = 1$) 
and $\bz_1=(1,1,\ldots,1)$. We denote by 
$\bz_1,\ldots,\bz_{\CD}$ the zeros of $\oU$. 
 
The reason for this choice is a strong presence of symmetries.
These symmetries guarantee that, for all $1\leq i,j\leq \CD$,
\begin{equation}\label{eq:s_1} 
   \mu(\oU,\bz_i) = \mu(\oU,\bz_j),
\end{equation}
and, consequently, that the value of $s_1$ is 
the same for all the zeros of $\oU$. Hence, 
\begin{eqnarray}\label{eq:U1}
   \E_{f\in\hd} K(f,\oU,\bz_1) &=& 
   \frac{1}{\CD}\sum_{j=1}^{\CD}\E_{f\in\hd}K(f,\oU,\bz_j)    
   \,=\,\E_{f\in\hd}\,\frac{1}{\CD}\sum_{j=1}^{\CD} K(f,\oU,\bz_j)\nonumber\\ 
   &=& \E_{f\in\hd}\,\frac{1}{\CD}\sum_{j=1}^{\CD} \Big(
     1+K(f,h_{s_1},\zeta^{(j)})\Big)\\
    &=& 1+H(\oU,\bz_1).\nonumber
\end{eqnarray}
That is, the average (w.r.t.~$f$) number of homotopy  
steps with initial system $\oU$ is the 
same no matter whether the zero of $\oU$ is taken at 
random or set to be $\bz_1$. Also, 
\begin{equation}\label{eq:U2}
    \mu^2(\oU,\bz_1) \le 2\, (n+1)^{D}
\end{equation}
(actually such bound holds for all zeros of $\oU$ but, again, 
this is not relevant for our argument). 
Both~\eqref{eq:s_1} and~\eqref{eq:U2} are proved 
in~\cite[Section~10.2]{Buergisser_Cucker-2010}. It follows 
from~\eqref{eq:U1},~\eqref{eq:homot3}, and~\eqref{eq:U2} 
that 
\begin{equation}\label{eq: complexityboundC2}
  \E_{f\in\hd} K(f,\oU,\bz_1) =\CO(D^3Nn^{D+1}).
\end{equation}
As above, for any fixed $\ep>0$ this bound is polynomial in $N$ 
provided $D\leq n^{\frac{1}{1+\ep}}$ and is 
bounded by $N^{\CO(\log\log N)}$ when    
$D\in[n^{\frac{1}{1+\ep}}, n]$.  
\eproof

\begin{lemma}\label{lem:t1}
With the notations of the proof of Theorem~\ref{th:near17} we have
\[
t_1=\frac{1}{\|f\|\sin\alpha \cot(s_1\alpha) -\|f\|\cos\alpha+1}
     \geq \frac{c'}{D^{3/2}\sqrt{N}\mu^2(g,\zeta)},
\]
$c'$ a constant.
\end{lemma}

\proof
The formula for $t_1$ is shown in~\cite[Prop.~5.2]{Buergisser_Cucker-2010}. 
For the bound, we have
\begin{eqnarray*}
  \|f\|\sin\a\cot(s_1\a)-\|f\|\cos\a+1 
&\leq & \|f\|\sin\a\, (s_1\a)^{-1} + \|f\|+1\\ 
&\leq & \sqrt{2N}\frac{1}{s_1} + \sqrt{2N}+1\\ 
&\leq & \sqrt{2N}\left(\frac{D^{3/2}\mu^2(g,\zeta)}{c} +1 
  +\frac{1}{\sqrt{2N}}\right)\\ 
&\leq & \frac{\sqrt{N}D^{3/2}\mu^2(g,\zeta)}{c'} 
\end{eqnarray*}
for an appropriately chosen $c'$.
\eproof

%%%%%%%%%%%%%

%%%%%%%%%%%

\end{document}